\documentclass[10pt,onecolumn,journal,papersize ]{article}
\usepackage[dvipdfmx]{graphicx}
\def\qed{\hfill \hbox{${\vcenter{\vbox{
   \hrule height 0.4pt\hbox{\vrule width 0.4pt height 6pt
   \kern5pt\vrule width 0.4pt}\hrule height 0.4pt}}}$}}
\newcommand{\maru}[1]{{\ooalign{\hfil#1\/\hfil\crcr
   \raise.167ex\hbox{\mathhexbox20D}}}}
\newcommand{\bysame}{%
   \leavevmode\hbox to 6em{\hrulefill}\,}

\begin{document}

\centerline{\large\bf On composite types of tunnel number two knots} 
\vskip 1mm 

\centerline{by} 
\vskip 2mm 

\centerline{\bf Kanji Morimoto} 
\vskip 5mm 

\centerline{Department of IS and Mathematics, Konan University} 
\vskip -1mm 

\centerline{Okamoto 8-9-1, Higashi-Nada, Kobe 658-8501, Japan}
\vskip -1mm 

\centerline{morimoto@konan-u.ac.jp} 
\vskip 5mm 

\begin{abstract}

Let $K$ be a tunnel number two knot. Then, by considering the $(g, b)$-decompositions, $K$ is one of (3, 0)-, (2, 1)-, (1, 2)- or (0, 3)-knots. In the present paper, we analyze the connected sum summands of composite tunnel number two knots and give a complete table of those summands from the point of view of $(g, b)$-decompositions. 
\vskip 2mm 

\vskip 2mm 

\end{abstract}

\vskip 10mm

{\bf 1. Introduction} 

Let $K$ be a knot in $S^3$. Then it is well known that $K$ can be uniquely decomposed into finitely many prime knots, which is due to Schubert ([11]) and is called the prime decomposition of $K$. Consider the tunnel number of $K$ denoted by $t(K)$, where the tunnel number is the minimal number of arcs properly embedded in the knot exterior $E(K)$ whose complementary space is homeomorphic to a handlebody. By the definition of the tunnel number, we have $t(K)=g(E(K))-1$, where $g(E(K))$ is the Heegaard genus of $E(K)$. In the present paper, we analyze the prime decomposition of tunnel number two knots and give a complete table from the point of view of $(g, b)$-decompositions defined as below. 

\hskip 4mm 
Let $g$ and $b$ be non-negative integers with $(g, b) \ne (0, 0)$. Then we say that a knot $K$ has a $(g, b)$-decomposition if there is a genus $g$ Heegaard splitting $(V_1, V_2)$ of $S^3$ such that $K$ intersects each handlebody in a $b$-string trivial arc system. In particular, if $b=0$ then we define a $(g, 0)$-decomposition as a genus $g$ Heegaard splitting of $S^3$ such that at least one of the two handlebodies contains $K$ as a central loop of a handle. If $g=0$, then a $(0, b)$-decomposition is the ordinary $b$-bridge decomposition. 

\hskip 4mm 
This decomposition is due to Doll ([1]) and is a generalization of the ordinary bridge decompositions due to Schubert ([12]). Then, by the definition and a little observation, we see that if a knot $K$ has a $(g, b)$-decomposition then $t(K) \le g+b-1$. This concept, $(g, b)$-decomposition of knots, plays very important role from the point of view of the tunnel numbers and the distance due to Hempel ([2]). For example, see [3], [4] or [8].  

\hskip 4mm 
Let $B$ be a 3-ball and $t_1 \cup t_2$ be two arcs properly embedded in $B$. Then $(B, t_1 \cup t_2)$ is called a 2-string tangle. We say that $(B, t_1 \cup t_2)$ is trivial if $t_1 \cup t_2$ is a 2-string trivial arc system in $B$, that $(B, t_1 \cup t_2)$ is free if $cl(B - N(t_1 \cup t_2))$ is a genus two handlebody, where $N(t_1 \cup t_2)$ is a regular neighborhood of $t_1 \cup t_2$ in $B$, that $(B, t_1 \cup t_2)$ is essential if $cl(\partial B - N(t_1 \cup t_2))$ is incompressible in  $cl(B - N(t_1 \cup t_2))$ and that $t_i \ (i=1,2)$ is unknotted if $(B, t_i)$ is a trivial ball pair. We say that a knot $K$ has a 2-string essential free tangle decomposition if $(S^3, K)$ is decomposed into a union of two 2-string essential free tangles. 

\hskip 4mm 
For a knot $K$ in $S^3$, we define the following two conditions $c(1)$ and $c(2)$. 
\vskip 3mm 

\hskip 4mm 
$c(1)$ : $(S^3, K)$ has a 2-string essential free tangle decomposition such that exactly one of the two tangles has an unknotted component. 

\hskip 4mm 
$c(2)$ : $(S^3, K)$ has a 2-string essential free tangle decomposition such that each tangle of the two tangles has an unknotted component. 
\vskip 3mm 

\hskip 4mm 
Under the above notations, for composite tunnel number two knots, we have shown the following : 
\vskip 3mm 

{\bf Theorem 1 ([5, 6])} \it Let $K$ be a composite tunnel number two knot, then one of the following holds. 

$(1)$ $K$ is the connected sum of a tunnel number one knot and a knot with a $(1, 1)$-decomposition, 

$(2)$ $K$ is the connected sum of a $2$-bridge knot and a knot with a $c(1)$- or a $c(2)$-condition. \rm  
\vskip 3mm 

\hskip 4mm  
For composite knots with $(2, 1)$-decompositions, we have shown the following : 
\vskip 3mm 

{\bf Theorem 2 ([7])} \it Let $K$ be a composite knot with a $(2, 1)$-decomposition, then one of the following holds. 

$(1)$ $K$ is the connected sum of a tunnel number one knot and a $2$-bridge knot,

$(2)$ $K$ is the connected sum of two knots with $(1, 1)$-decompositions, 

$(3)$ $K$ is the connected sum of a $2$-bridge knot and a knot with a $c(2)$-condition. \rm 
\vskip 3mm 

\hskip 4mm 
For 3-bridge knots, by the additivity of bridge indices due to Schubert, we have the following : 
\vskip 3mm 

{\bf Theorem 3 ([12])} \it Let $K$ be a composite $3$-bridge knot, then $K$ is the connected sum of two $2$-bridge knots.   \rm 
\vskip 3mm 

\hskip 4mm 
In the present paper, for knots with $(1, 2)$-decompositions, we will show the following : 
\vskip 3mm 

{\bf Theorem 4} \it Let $K$ be a composite knot with a $(1, 2)$-decomposition, then $K$ is the connected sum of a knot with a $(1, 1)$-decomposition and a $2$-bridge knot. \rm 
\vskip 3mm 

\hskip 4mm 
By the way, to state the above results more precisely, we need to define the term $(g, b)$-knots. To do this, by a little observation, we have : 
\vskip 3mm 

{\bf Fact 1} \it $(1)$ If a knot $K$ has a $(g-1, b+1)$-decomposition, then $K$ has a $(g, b)$-decomposition. 

$(2)$ If a knot $K$ has a $(g, b-1)$-decomposition, then $K$ has a $(g, b)$-decomposition. \rm 
\vskip 3mm 

\hskip 4mm 
By the above fact, in the present paper, we say that a knot $K$ is a $(g, b)$-knot if $K$ has a $(g, b)$-decomposition but has neither $(g-1, b+1)$-decomposition nor $(g, b-1)$-decomposition. In addition, we say that a knot $K$ is a $c(i)$-knot if $K$ has a $c(i)$-condition $(i=1,2)$. Then the above Theorems 1, 2, 3 and 4 are rewritten as follows. We note that, since $(3, 0)$-knots are tunnel number two knots which have no $(2, 1)$-decompositions, it is needed to delete the $c(2)$-condition in Theorem 1. We further note that there is no knot which has both conditions $c(1)$ and $c(2)$ because of the unique 2-string essential free decomposition theorem due to Ozawa ([10]).  
\vskip 3mm 

{\bf Theorem 1} \it Let $K$ be a composite $(3, 0)$-knot, then one of the following holds. 

$(1)$ $K$ is the connected sum of a $(2, 0)$-knot and a $(1, 1)$-knot, 

$(2)$ $K$ is the connected sum of a $(0, 2)$-knot and a $c(1)$-knot. \rm 
\vskip 3mm 

{\bf Theorem 2} \it Let $K$ be a composite $(2, 1)$-knot, then one of the following holds. 

$(1)$ $K$ is the connected sum of a $(2, 0)$-knot and a $(0, 2)$-knot,

$(2)$ $K$ is the connected sum of two $(1, 1)$-knots, 

$(3)$ $K$ is the connected sum of a $(0, 2)$-knot and a $c(2)$-knot. \rm 
\vskip 3mm 

{\bf Theorem 3} \it Let $K$ be a composite $(0, 3)$-knot, then $K$ is the connected sum of two $(0, 2)$-knots. \rm 
\vskip 3mm 

{\bf Theorem 4} \it Let $K$ be a composite $(1, 2)$-knot, then $K$ is the connected sum of a $(1, 1)$-knot and a $(0, 2)$-knot. \rm 
\vskip 3mm 

\hskip 4mm
Then, by summarizing the above results, we have the following table of composite tunnel number two knots from the point of view of $(g, b)$-decompositions.   

\begin{center}
\begin{tabular}
{@{\vrule width0.8pt~}c|l@{~\vrule width0.8pt}}
\noalign{\hrule height0.8pt}
$(g, b)$ & composite types \\
\noalign{\hrule height0.8pt}
$(3, 0)$ & $(2, 0) \# (1, 1)$ \ or \ $(0, 2) \# c(1)$ \\
\hline 
$(2, 1)$ & $(2, 0) \# (0, 2)$, \ $(1, 1) \# (1, 1)$ \ or \ $(0, 2) \# c(2)$ \\
\hline 
$(1, 2)$ & $(1, 1) \# (0, 2)$ \\
\hline 
$(0, 3)$ & $(0, 2) \# (0, 2)$ \\
\noalign{\hrule height0.8pt}
\end{tabular}
\vskip 3mm 
Table 1: Composite types of tunnel number two knots
\end{center}
\vskip 3mm 

{\bf Remark 1} Concerning the $c(1)$- and $c(2)$-conditions, we have the following : 

(1) If a knot $K$ has $c(1)$-condition then $K$ is a prime $(3, 0)$-knot. Hence $(0, 2) \# c(1)$ is included in $(0, 2) \# (3, 0)$ and this is the tunnel number degeneration $``2+1=2"$. 

(2) If a knot $K$ has  $c(2)$-condition then $K$ is a prime $(2, 1)$-knot. Hence $(0, 2) \# c(2)$ is included in $(0, 2) \# (2, 1)$ and this is also the tunnel number degeneration $``2+1=2"$. 
\vskip 3mm 

{\bf Remark 2} By Fact 1, the family of knots with $(g, b)$-decompositions contains the family of knots with $(g-1, b+1)$-decompositions and the family of knots with $(g, b-1)$-decompositions. Hence we have the $(g, b)$-diagram as illustrated in Figure 1. 

\begin{figure}[htbp]
\centerline{\includegraphics[width=7.5cm]{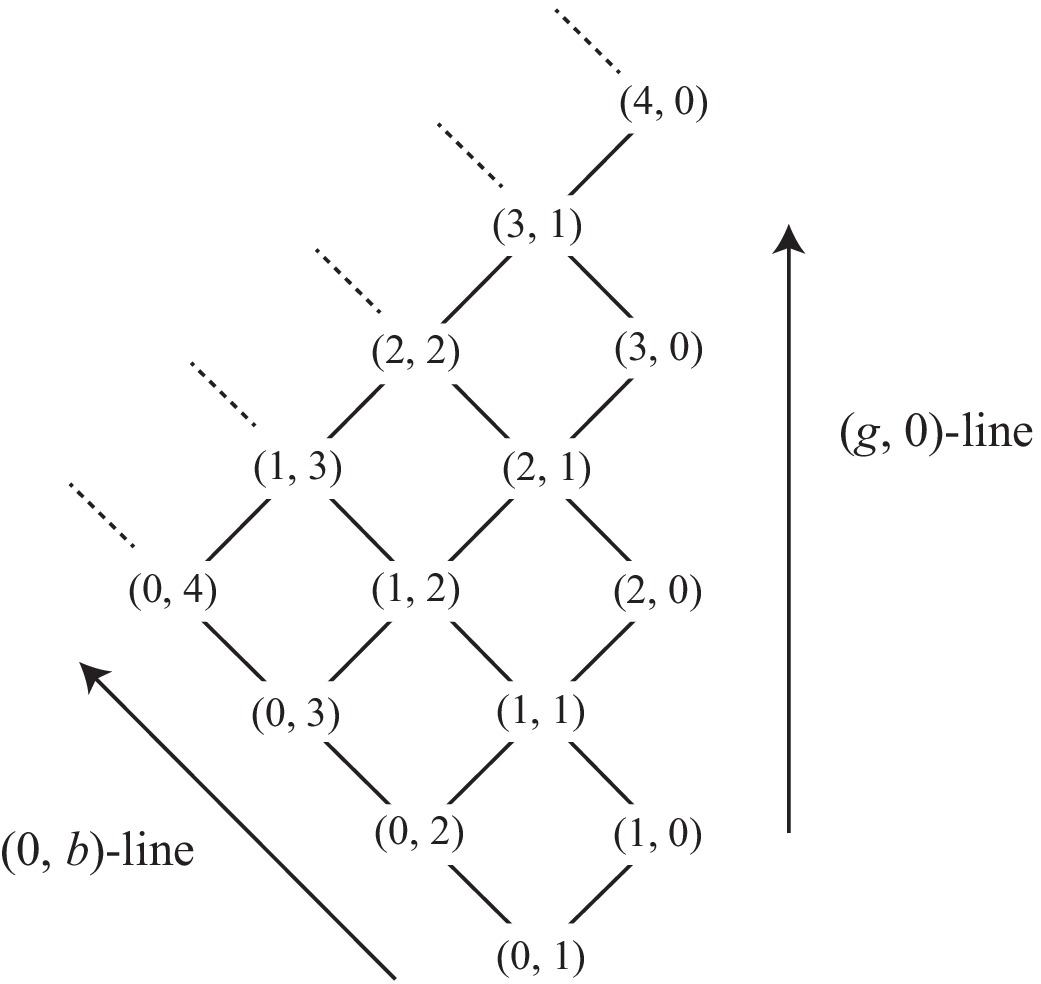}}
\vskip 3mm 

\centerline{Figure 1: The $(g, b)$-diagram}
\end{figure}
\vskip 10mm

{\bf 2. Proof of Theorem 4}

Let $K$ be a composite knot in $S^3$ with the decomposing 2-sphere $S$. Suppose $K$ has a $(1, 2)$-decomposition. Then there is a genus one Heegaard splitting $(V_1, V_2)$ of $S^3$ such that $K$ intersects $V_i \ (i=, 2)$ in 2-string trivial arc system, where $V_i$ is a solid torus.   

\hskip 4mm 
Put $V_i \cap K = \gamma_i^1 \cup \gamma_i^2$, and $S_i = V_i \cap S$. Then, by taking a spine of $V_1$, we may assume that $S_1$ consists of disks not intersecting $\gamma_1^1 \cup \gamma_1^2$ and $S_2$ is a planar surface properly embedded in $V_2$ intersecting $\gamma_2^1 \cup \gamma_2^2$ in two points. In addition we may assume that the number of the components of $S_1$ is minimal among all decomposing 2-spheres of $K$. Then we have : 
\vskip 3mm 

{\bf Lemma 2.1} \it Each component of $S_1$ is one of the following three types as in Figure $2$: 

$(1)$ a separating disk which cuts off a $3$-ball containing one of $\gamma_1^1 \cup \gamma_1^2$,  

$(2)$ a separating disk which cuts off a $3$-ball containing both of $\gamma_1^1 \cup \gamma_1^2$, 

$(3)$ a non-separating disk. \rm 

{\bf Proof.} Let $D$ be a component of $S_1$. Suppose $D$ is a separating disk. Then $D$ divides $V_1$ into a 3-ball and a solid torus. If the 3-ball contains no component of $\gamma_1^1 \cup \gamma_1^2$, then we can reduce the number of the components of $S_1$ and this contradicts the minimality. Thus the 3-ball contains at least one component of $\gamma_1^1 \cup \gamma_1^2$ and this completes the proof of the lemma. \qed  

\begin{figure}[htbp]
\centerline{\includegraphics[width=13cm]{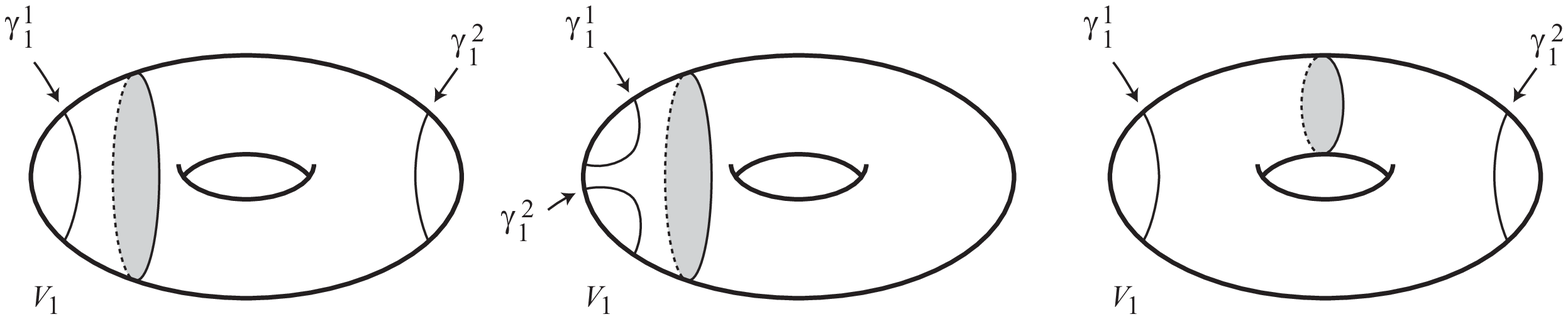}}

\centerline{(1) \hskip 4cm (2) \hskip 4cm (3)}

\centerline{Figure 2: Disks in $V_1$}
\end{figure}

\hskip 4mm 
Next, let $E_1$ and $E_2$ the disks in $V_2$ for the triviality of $\gamma_2^1$ and $\gamma_2^2$ respectiverly, and $E_3$ and $E_4$ the two non-separating disks in $V_2$ such that $E_3 \cup E_4$ divides $V_2$ into two 3-balls each of which contains one of $E_1 \cup E_2$ as in Figure 3. 
 
\begin{figure}[htbp]
\centerline{\includegraphics[width=5cm]{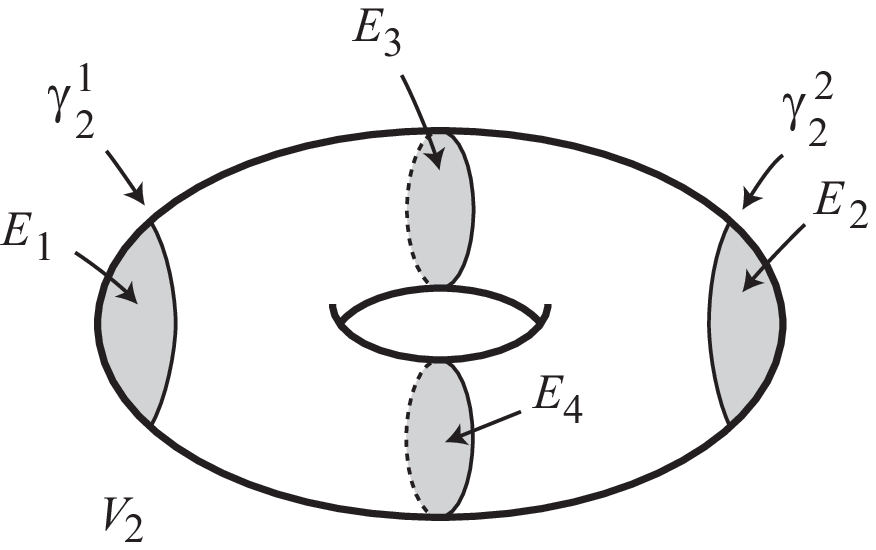}}

\centerline{Figure 3: $E_1 \cup E_2 \cup E_3 \cup E_4$ in $V_2$}
\end{figure}
 
\hskip 4mm 
Put $E = E_1 \cup E_2 \cup E_3 \cup E_4$, then by standard cut and paste operations, we may assume that each component of $S_2 \cap E$ is an arc properly embedded in $E$. We say that an arc $\alpha$ properly embedded in $S_2$ is $\gamma$-essential if $\alpha$ is essential in $S_2 - (\gamma_2^1 \cup \gamma_2^2)$. Suppose there is an arc $\alpha$ in $S_2 \cap E$ which is $\gamma$-essential in $S_2$. Let $\Delta$ be the disk in $E$ such that $\partial \Delta$ is the union of $\alpha$ and a subarc of $\partial E - (\gamma_2^1 \cup \gamma_2^2)$. We may assume that $\Delta \cap S_2 = \alpha$ by changing the disks $E$ if necessary. Then we can perform a boundary compression of $S_2$ at $\alpha$ along $\Delta$ from $V_2$ to $V_1$, and we get a band, say $b$, in $V_1$. If $b$ connects two different disks in $S_1$, then we can reduce the number of the components of $S_1$ and this contradicts the minimality. Thus $b$ connects a single disk, and the union of the band and the disk is an annulus in $V_1$. Then we have : 
\vskip 3mm 

{\bf Lemma 2.2} \it The annulus is one of the following three types as in Figure $4$ : 

$(1)$ the union of a separating disk of type $(1)$ in Lemma $2.1$ and a band which is contained in the solid torus component and winds around the longitude exactly once, 

$(2)$ the union of a separating disk of type $(2)$ in Lemma $2.1$ and a band which is contained in the 3-ball component such that the compressing disk of the annulus in $\partial V_1$ intersects $\gamma_1^1 \cup \gamma_1^2$ in two points.  

$(3)$ the union of a non-separating disk and a band such that the compressing disk of the annulus in $\partial V_1$ intersects $\gamma_1^1 \cup \gamma_1^2$ in two points, \rm 

\begin{figure}[htbp]
\centerline{\includegraphics[width=13cm]{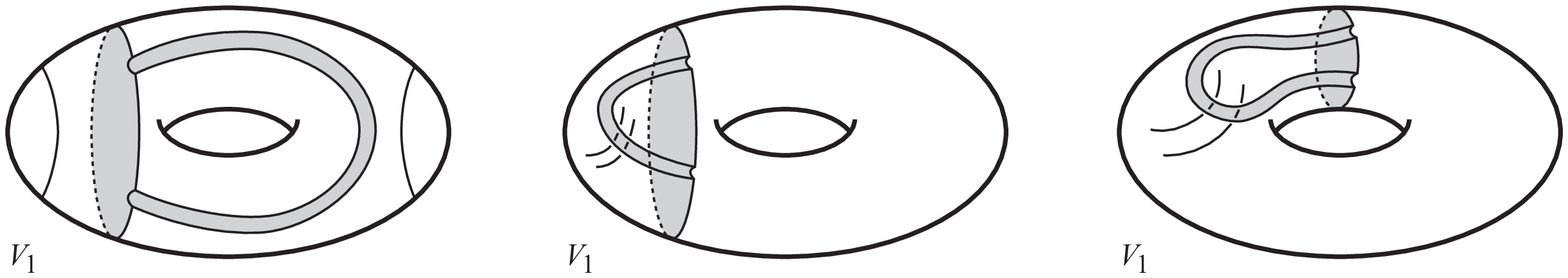}}

\centerline{(1) \hskip 4cm (2) \hskip 4cm (3)}

\centerline{Figure 4: Annuli in $V_1$}
\end{figure}

{\bf Proof.} Let $D$ be a disk component of $S_1$ which is connected to the band $b$, then $D \cup b$, say $A$, is an annulus properly embedded in $V_1$. Suppose $D$ is a separating disk of type (1). Then, since $D$ divides $V_1$ into a 3-ball and a solid torus, $b$ is contained in the 3-ball or in the solid torus. If $b$ is contained in the 3-ball, then $A$ is compressible in $V_1$ and by standard cut and paste operation, we can take another decomposing 2-sphere which intersects $V_1$ in fewer disk components. This contradicts the minimality. Thus $b$ is contained in the solid torus component. Then by the same reason as above, $b$ winds around the solid torus at least once in the longitudinal direction. However, if the band winds more than once, then, since we can regard each disk of $S-A$ is a 2-handle for $V_1$, we have the lens space not $S^3$. This is a contradiction. Thus $b$ winds around $V_1$ exactly once and $A$ is an annulus of type (1). 

\hskip 4mm 
Next suppose $D$ is a separating disk of type (2). Then, since $D$ divides $V_1$ into a 3-ball and a solid torus, $b$ is contained in the 3-ball or in the solid torus. If $b$ is contained in the solid torus, then by the same reason as above $b$ winds around the solid torus exactly once. However, in this case, $A$ can be pushed out into $V_2$ and we can reduce the number of the components of $S_1$. Thus $b$ is contained in the 3-ball component. Then each component of $\partial A$ bounds a disk in $\partial V_1$. If at least one disk intersects $\gamma_1^1 \cup \gamma_1^2$ in 0 or 1 point, then we can take another decomposing 2-sphere which intersects $V_1$ in fewer disk components. This contradicts the minimality and shows that $A$ is an annulus of type (2). 

\hskip 4mm 
Finally suppose $D$ is a non-separating disk. Then, since $D$ cuts open $V_1$ into a 3-ball, $A$ is compressible in $V_1$. Let $\Delta$ be a compressing disk for $A$ in $\partial V_1$. If $\Delta \cap (\gamma_1^1 \cup \gamma_1^2) =$ 0 or 1 point, then by cut and paste operation we can take another decomposing 2-sphere which intersects $V_1$ in fewer disk components. This contradicts the minimality. Thus $\Delta \cap (\gamma_1^1 \cup \gamma_1^2) =$ 2, 3 or 4 points. However, if it is 3 or 4 points, then by taking another compressing disk for $A$, i.e., a meridian disk of $V_1$, we can reduce the number of the components of $S_1$ and have a contradiction similarly. Thus $A$ is an annulus of type (3), and this completes the proof. \qed 
\vskip 3mm 

\hskip 4mm 
Let $n$ be the number of the components of $S_1$, then : 
\vskip 3mm 

{\bf Lemma 2.3} \it We have $n=1$. Hence $S_1$ is a single disk not intersecting $\gamma_1^1 \cup \gamma_1^2$ and $S_2$ is a single disk intersecting $\gamma_2^1 \cup \gamma_2^2$ in two points. \rm 

{\bf Proof.} Suppose $n>1$. Then, since $S_2$ has a non-trivial fundamental group, $S_2 \cap E$ has a $\gamma$-essential arc. Thus we can perform a boundary compression and get a band $b_1$ in $V_1$. Then by the minimality of $S_1$, $b_1$ connects a single component of $S_1$ and an annulus is produced in $V_1$. Continue this procedure. Then, at some $k$th stage, we have that $b_1, b_2, \cdots, b_k$ are the bands each of which connects a single component of $S_1$ and $b_{k+1}$ connects two different components of $S_1$, where to avoid the confusion of notations we use the same notations of $S_1$ and $S_2$ even after the boundary compressions. We note that there are no two bands $b_i$ and $b_j$ which connect the same disk because at each stage core arcs of those bands are $\gamma$-essential in $S_2$. Then, at this stage, we have $k$ annuli $A_1, A_2, \cdots, A_k$ in $V_1$. 

\hskip 4mm 
Suppose $k<n$. Then, since there remains a disk component of $S_1$, we have the following two cases : 

\hskip 4mm 
(i) $b_{k+1}$ connects a disk and an annulus $A_i$, 

\hskip 4mm 
(ii) $b_{k+1}$ connects two annuli $A_i, A_j (i<j)$. 

\hskip 4mm 
In case (i), we can use the inverse operation of boundary compression introduced by Ochiai in [9], and can reduce the number of the components of $S_1$. This is a contradiction. In case (ii), we have two subcases, (ii-a) : $A_1, A_2, \cdots, A_k$ are all mutually parallel annuli of type (1) in Lemma 2.2, (ii-b) : $A_1, A_2, \cdots, A_k$ are of type (2) or of type (3) in Lemma 2.2. Then, in case (ii-a), $b_{k+1}$ does not run over the band $b_j$. In case (ii-b), by the existence of the disk of type (3) in Lemma 2.1, we see that $b_{k+1}$ does not run over the band $b_j$. Then, in both cases, we can pull back the bands $b_j, \cdots, b_k$ leaving $b_{k+1}$in $V_1$. This means that case (ii) is reduced to case (i) and we have a contradiction. Thus we have $k=n$. 

\hskip 4mm 
By the above arguments, we can put $S_1 = A_1 \cup A_2 \cup \cdots \cup A_n$ and $S_2 = D_1^* \cup D_2^* \cup B_1 \cup B_2 \cup \cdots \cup B_{n-1}$, wher $A_i \ (i=1, 2, \cdots, n)$ is an annulus in $V_1$, $B_i \ (i=1, 2, \cdots, n-1)$ is an annulus in $V_2$ and $D_i^* \ (i=1, 2)$ is a disk in $V_2$ intersecting $\gamma_2^1 \cup \gamma_2^2$ in a single point.   

\hskip 4mm 
Suppose $D_1^*$ is a non-separating disk and $\partial D_1^*$ is identified with a component of $\partial A_i$ for some $i$. Then, since $V_1 \cup V_2 = S^3$, $A_i$ is of type (1) in Lemma 2.2. Then, by Lemma 2.2, those components of $S_1$ are all mutually parallel annuli of type (1). This means that each component of $\partial S_2$ is a meridian of $V_2$. However, each annulus component of $S_2$ has a  boundary component which is not a meridian by Lemma 2.2. This is a contradiction, and shows that both of $D_1^*$ and $D_2^*$ are separating disks. 

\hskip 4mm 
Let $X$ be a 3-ball in $V_2$ cut off by $D_1^*$. Then we may assume that $X \cap S_2 = D_1^*$, and since $D_i^* \cap (\gamma_2^1 \cup \gamma_2^2) \ (i=1,2)$ is a single point, $(X \cap \partial V_2) \cap (\gamma_2^1 \cup \gamma_2^2)$ is a single point or three points. However, there is no annulus component of $S_1$ whose boundary component bounds a disk in $\partial V_1$ intersecting $\gamma_1^1 \cup \gamma_1^2$ in a single point or in three points. This contradiction is due to the hypothesis $n>1$, and completes the proof. \qed  
\vskip 3mm 

{\bf Proof of Theorem 4.} By Lemma 2.3, $S_1$ is a single separating disk not intersecting $\gamma_1^1 \cup \gamma_1^2$, and $S_2$ is a single separating disk intersecting $\gamma_2^1 \cup \gamma_2^2$ in two points. We may assume that $S_2 \cap E$ consists of arcs properly embedded in $E$ and that the number of the components of $S_2 \cap E$ is minimal among all decomposing 2-spheres intersecting $V_i \ (i=1,2)$ in a single disk. We note that $S_2 \cap E$ contains two arcs which meets a point of $S_2 \cap (\gamma_2^1 \cup \gamma_2^2)$. Then we have the following two cases : 

\hskip 4mm 
(i) $S_2 \cap E$ consists of exactly two arcs each of which meets a point of $S_2 \cap (\gamma_2^1 \cup \gamma_2^2)$, 

\hskip 4mm 
(ii) $S_2 \cap E$ contains a $\gamma$-essential arc properly embedded in $S_2$. 

\hskip 4mm    
Suppose we are in case (i). Then we have an arc component $\alpha$ of $S_2 \cap E_1$ which intersects $\gamma_2^1$ in a single point. Let $\Delta$ be one of the two disks in $E_1$ cut off by $\alpha$. Then we may assume $\Delta \cap S = \alpha$, and can isotope $S$ along the disk $\Delta$ so that $S \cap V_i \ (i=1,2)$ is a single separating disk intersecting $\gamma_i^1 \cup \gamma_i^2$ in a single point as in Figure 5. 

\begin{figure}[htbp]
\centerline{\includegraphics[width=9cm]{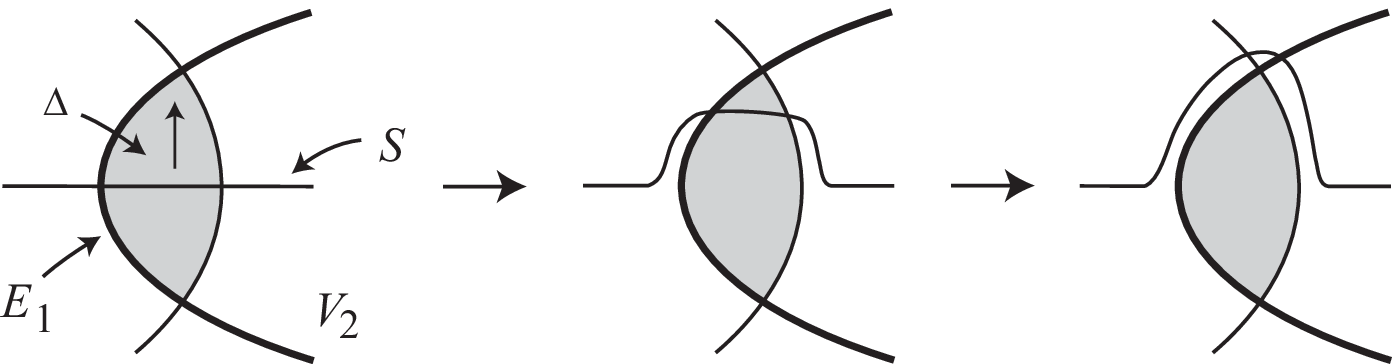}}

\centerline{Figure 5: The isotopy along $\Delta$}
\end{figure}

\hskip 4mm
Let $X_i \ (i=1,2)$ be the 3-ball in $V_i$ cut off by $S_i$. Since, $S_i \cap (\gamma_i^1 \cup \gamma_i^2)$ is a single point, $(X_i \cap \partial V_i) \cap (\gamma_i^1 \cup \gamma_i^2)$ is one or three points. If it is one point, then the tangle $(X_i, X_i \cap (\gamma_i^1 \cup \gamma_i^2))$ is a trivial tangle. In this case, $(X_1 \cup X_2, (X_1 \cup X_2) \cap K)$ is a 1-string trivial tangle and this means that $S$ bounds a trivial ball pair. This is a contradiction because $S$ is a decomposing 2-sphere of the non-trivial connected sum of $K$. Thus we have that $(X_i \cap \partial V_i) \cap (\gamma_i^1 \cup \gamma_i^2)$ consists of three points as in Figure 6. 

\begin{figure}[htbp]
\centerline{\includegraphics[width=10cm]{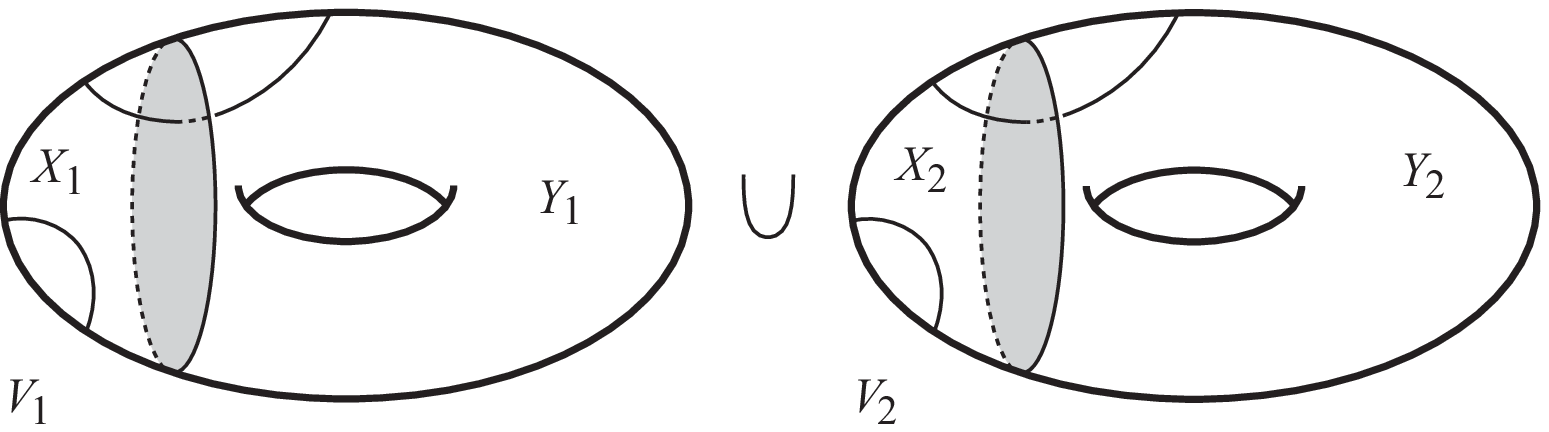}}

\centerline{Figure 6: The intersection of $S$ and $V_1 \cup V_2$ in case (i)}
\end{figure}

\hskip 4mm 
Let $X_i \ (i=1,2)$ be the 3-ball in $V_i$ cut off by $S_i$ and $Y_i$ the solid torus in $V_i$ cut off by $S_i$. Then $(X_i, X_i \cap (\gamma_i^1 \cup \gamma_i^2))$ is a 2-string trivial tangle and $(Y_i, Y_i \cap (\gamma_i^1 \cup \gamma_i^2))$ is a solid torus with a single trivial arc. Hence $X_1 \cup X_2$ extends to a 2-bridge decomposition of a knot and $Y_1 \cup Y_2$ extends to a knot with a (1, 1)-decomposition, i.e., $K$ is the connected sum of a (0, 2)-knot and a (1, 1)-knot. 

\hskip 4mm 
Next, suppose we are in case (ii). Let $\alpha$ be a $\gamma$-essential arc properly embedded in $S_2$. Then we may assume that $\alpha$ is outermost in $E$. Perform a boundary compression of $S_2$ at $\alpha$. Then, since $\alpha$ is an arc in $S_2$ which splits the two points $S_2 \cap (\gamma_2^1 \cup \gamma_2^2)$, we can put $S_1 = A$ and $S_2 = D_1^* \cup D_2^*$, where $A$ is an annulus of type (1) in Lemma 2.2 and $D_i^* \ (i=1,2)$ is a meridian disk of $V_2$ intersecting $\gamma_2^1 \cup \gamma_2^2$ in a single point. We note that if $A$ is an annulus of type (2) in Lemma 2.2, then $D_i^*$ is a separating disk intersecting $\gamma_i^1 \cup \gamma_i^2$ in two points and this is a contradiction. 

\hskip 4mm 
Let $X_i \cup Y_i \ (i=1,2)$ be the two components cut off by $S_i$, where $X_1$ and $X_2$ are identified and $Y_1$ and $Y_2$ are identified as in Figure 7. Then, by adding a 2-handle to $X_1$ along the annulus $A$, $(X_1, \gamma_1^1)$ extends to a 2-string trivial tangle, where a core arc of the 2-handle is regarded as a string. In addition $(X_2, X_2 \cap (\gamma_2^1 \cup \gamma_2^2))$ is also a 2-string trivial tangle. 

\begin{figure}[htbp]
\centerline{\includegraphics[width=10cm]{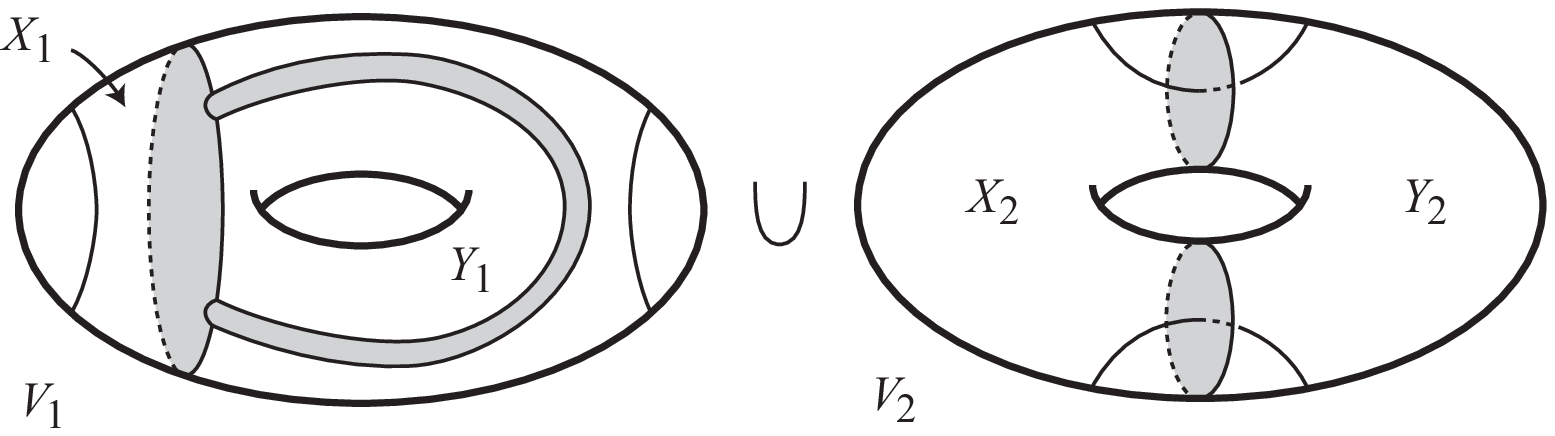}}

\centerline{Figure 7: The intersection of $S$ and $V_1 \cup V_2$ in case (ii)}
\end{figure}

\hskip 4mm 
Hence $X_1 \cup X_2$ extends to a 2-bridge decomposition of a knot. On the other hand, $(Y_1, \gamma_1^2)$ is a solid torus with a trivial arc, and by adding a 1-handle to $Y_2$ along the two disks $D_1^* \cup D_2^*$, $(Y_2, Y_2 \cap (\gamma_2^1 \cup \gamma_2^2))$ extends to a solid torus with a trivial arc, where a core arc of the 1-handle is regarded as a part of the trivial arc. Hence $Y_1 \cup Y_2$ extends to a $(1, 1)$-decomposition of a knot. This shows that $K$ is the connected sum of a (0, 2)-knot and a (1, 1)-knot, and completes the proof of Theorem 4. \qed    
\vskip 10mm

{\bf\large References}
\vskip 3mm 

\ [1] \ H. Doll, \ A generalized bridge number for links in 3-manifolds, Math. Ann., {\bf 294} 

\hskip 7mm (1992) 701--717

\ [2] \ J. Hempel, \ 3–manifolds as viewed from the curve complex. Topology, {\bf 40} 

\hskip 7mm (2001) 631-–657.

\ [3] \ J. Johnson and A. Thompson, \ On tunnel number one knots that are not $(1, n)$, 

\hskip 7mm J. Knot Theory Ramifications, {\bf 20} (2011) 609--615. 

\ [4] \ T. Li and R. Qiu, \ On the degeneration of tunnel numbers under connected 

\hskip 7mm sum, arXiv:math/1310.5054.

\ [5] \ K. Morimoto, \ On the additivity of tunnel number of knots, Topology Appl., 

\hskip 7mm {\bf 53} (1993) 37--66.

\ [6] \ K. Morimoto, Charaterization of tunnel number two knots which have the prop- 

\hskip 7mm erty  `` $2 + 1 = 2$ '', Topology Appl., {\bf 64} (1995) 165-176.

\ [7] \ K. Morimoto, Charaterization of composite knots with 1-bridge genus two,  

\hskip 7mm J. Knot Theory Ramifications, {\bf 10} (2001) 823--840.

\ [8] \ Y. Minsky, Y. Moriah and S. Schleimer, \ High distance knots, arXiv:math/0607265. 

\ [9] \ M. Ochiai, On Haken's theorem and its extension, Osaka J. Math. {\bf 20} (1983) 

\hskip 7mm 461--480.

[10] \ M. Ozawa, On uniqueness of essential tangle decompositions of knots with free 

\hskip 7mm tangle decompositions, Proc. Appl. Math. Workshop {\bf 8} ed. G. T. Jin and 

\hskip 7mm K. H. Ko, KAIST, Taejon (1998) 227--232. 

[11] \ H. Schubert, Die eindeutige Zerlegbarkeit eines Knoten in Primknoten, Sitzungs-

\hskip 7mm ber. Akad. Wiss. Heidelberg, math.-nat. KI. 3 (1949) Abh: 57--104.

[12] \ H. Schubert, $\ddot U$ber eine numerische Knoteninvariante, Math. Z. {\bf 61} (1954) 245-

\hskip 7mm 288.

\end{document}